\documentclass{article}

\usepackage{amssymb}
\usepackage{amsmath}
\usepackage{amsfonts}
\usepackage{graphicx}
\usepackage{hyperref}%
\providecommand{\U}[1]{\protect\rule{.1in}{.1in}}
\RequirePackage{doi}
\newtheorem{theorem}{Theorem}[section]
\begin{document}

\title{Isometric flows of $G_2$-structures}

\author{Sergey Grigorian\\
University of Texas Rio Grande Valley\\
1201 W. University Drive\\
Edinburg, TX 78541\\
USA}
%



%

\maketitle

\begin{abstract}
We survey recent progress in the study of flows of isometric $G_2$-structures on 7-dimensional manifolds, that is, flows that preserve the metric, while modifying the $G_2$-structure. In particular, heat flows of isometric $G_2$-structures have been recently studied from several different perspectives, in particular in terms of $3$-forms, octonions, vector fields, and geometric structures. We will give an overview of each approach, the results obtained, and compare the different perspectives. 
\end{abstract}   


%

\section{Introduction}

\label{secIntro}One of the most challenging problems
in differential geometry is the question of existence conditions for
torsion-free $G_{2}$-structures on smooth $7$-dimensional manifolds. Such
$G_{2}$-structures are precisely the ones that correspond to metrics with
holonomy contained in $G_{2}$. One approach that has been pioneered by Robert
Bryant \cite{bryant-2003} is to considered heat-like flows of $G_{2}%
$-structures with the hope that under certain conditions they may converge to
a torsion-free $G_{2}$-structure. A difficulty that is encountered in such an
approach is that in general, deformations of a $G_{2}$-structure also affect
the corresponding metric, and so any heat equation for the $G_{2}$-structure
becomes nonlinear. This is not unlike the situation for the Ricci flow, where
the underlying geometry changes along the flow, however in the $G_{2}$ case,
we have two separate but closely related objects, the $G_{2}$-structure and
the metric, both of which vary along the flow. Given a Riemannian metric on a
$7$-manifold that admits $G_{2}$-structures, there is a family of $G_{2}%
$-structures that correspond to it, so a possible approach could be to
separate as much as possible the deformations of the metric from the
deformations of $G_{2}$-structures that preserve the metric. Indeed, as was
shown by Karigiannis \cite{karigiannis-2005-57}, given a decomposition of
$3$-forms according to representations of $G_{2}$, the deformations of the
$G_{2}$-structure $3$-form that preserve the metric are precisely the ones that lie in the
$7$-dimensional representation $\Lambda_{7}^{3}$. Bryant's original Laplacian
flow of closed $G_{2}$-structures has no component in $\Lambda_{7}^{3}$
\cite{bryant-2003}, and as such is transverse to directions that preserve the
metric. This allowed for more tractable analytic properties. In contrast, a
similar flow for co-closed $G_{2}$-structures that was proposed in
\cite{KarigiannisMcKayTsui} does have a component in $\Lambda_{7}^{3},$ which,
as shown in \cite{GrigorianCoflow}, causes non-parabolicity of the flow. This
suggests that the freedom of $G_{2}$-structures to move in directions that
preserve the metric is some kind of degeneracy and thus suitable gauge-fixing
conditions within the metric class are needed to address it.

These considerations show that it is necessary to have a clearer picture of
$G_{2}$-structures within a fixed metric class. In \cite{bryant-2003}, Bryant
observed that such $G_{2}$-structures are parametrized by sections of an
$\mathbb{R}P^{7}$-bundle, or more concretely, by pairs $\left(  a,\alpha
\right)  $ where $a$ is a real-valued function and $\alpha$ is a vector field
such that $a^{2}+\left\vert \alpha\right\vert ^{2}=1,$ and $\pm\left(
a,\alpha\right)  $ define the same $G_{2}$-structure. If $\varphi$ is a fixed
$G_{2}$-structure, then any other $G_{2}$-structure $\sigma_{\left(
a,\alpha\right)  }\left(  \varphi\right)  $ within the same metric class is
given by:
\begin{equation}
\sigma_{\left(  a,\alpha\right)  }\left(  \varphi\right)  =\left(
a^{2}-\left\vert \alpha\right\vert ^{2}\right)  \varphi-2a\alpha\lrcorner
\psi+2\alpha\wedge\left(  \alpha\lrcorner\varphi\right)  ,\label{phisameg1}%
\end{equation}
where $\psi=\ast\varphi$.

Given that the group $G_{2}$ may be defined as the automorphism group of the
octonions, a $G_{2}$-structure defines an octonion structure on the manifold,
and in \cite{GrigorianOctobundle}, this observation was used to interpret the
above pair $\left(  a,\alpha\right)  $ as a unit octonion $V$, and then
(\ref{phisameg1}) is just the $3$-form that corresponds to a modified octonion
product defined by $V.$ Thus, a flow of isometric $G_{2}$-structures can be
interpreted as a flow of the unit octonion section $V.$ In particular, a
natural heat flow of isometric $G_{2}$-structures was introduced in
\cite{GrigorianOctobundle}. Given an octonionic covariant derivative $D$,
constructed from the Levi-Civita connection and the torsion of the initial
$G_{2}$-structure $\varphi$, the heat flow of isometric $G_{2}$-structures is
then the semilinear, parabolic equation
\begin{equation}
\frac{\partial V}{\partial t}=\Delta_{D}V+\left\vert DV\right\vert ^{2}V
\label{octoflow}%
\end{equation}
with some initial condition $V\left(0\right)  =V_{0}$ and where $\Delta
_{D}=-D^{\ast}D$ is the Laplacian operator corresponding to $D$. This was
obtained as the negative gradient flow of an energy functional with respect to
$D$. The critical points of the flow (\ref{octoflow}) correspond to $G_{2}%
$-structures for which the torsion tensor is divergence-free, i.e. satisfies
$\operatorname{div}T=0,$ where divergence is taken with respect to the
Levi-Civita connection. This is significant for several reasons. The
divergence of torsion is precisely the term that causes the non-parabolicity
of the Laplacian flow of co-closed $G_{2}$-structures from \cite{KarigiannisMcKayTsui}  as mentioned above, and $\operatorname{div}%
T=0$ for closed $G_{2}$-structures. Thus, closed $G_{2}$-structures are
automatically critical points of (\ref{octoflow}). Secondly, $T$ has been
interpreted in \cite{GrigorianOctobundle} as an imaginary octonion-valued
$1$-form, which is added to the Levi-Civita connection to obtain the
octonionic covariant derivative $D,$ hence the condition $\operatorname{div}%
T=0$ is precisely analogous to the Coulomb gauge condition in gauge theory.
This analogy makes this condition a reasonable candidate for a gauge-fixing
condition within a fixed metric class.

Soon after the introduction of the flow (\ref{octoflow}) in
\cite{GrigorianOctobundle}, it was further studied from different perspectives
by several authors: Bagaglini in \cite{Bagaglini2}; Dwivedi,
Gianniotis, and Karigiannis in \cite{DGKisoflow}; the author in
\cite{GrigorianIsoflow}; Loubeau and S\'{a} Earp in \cite{SaEarpLoubeau}.

Equivalently to the flow of octonions (\ref{octoflow}), one can consider
directly the evolution of the $3$-form $\varphi$ via the equation
\begin{equation}
\frac{\partial\varphi}{\partial t}=2\left(  \operatorname{div}T\right)
\lrcorner\psi\label{octoflow2}%
\end{equation}
where $T$ is the torsion tensor that corresponds to the $G_{2}$-structure
$3$-form at time $t.$ This is the way the flow was formulated in
\cite{Bagaglini2} and in \cite{DGKisoflow} (although here we are
following \cite{GrigorianOctobundle, GrigorianIsoflow} and added a factor of 2
in (\ref{octoflow2}). In \cite{SaEarpLoubeau}, a more general approach is
taken and a harmonic heat flow of geometric structures is considered. In the
case of $G_{2}$-structures, it is shown to reduce to (\ref{octoflow2}). In
this survey we will review the above approaches to the flow of isometric
$G_{2}$-structures and outline the key analytic results.
\subsection*{Acknowledgements}
This work was supported by the National Science Foundation [DMS-1811754].

\section{Isometric $G_{2}$-structures}

\label{secIso}A $G_{2}$-structure on a $7$-manifold is
defined by a smooth \emph{positive} $3$-form $\varphi$
\cite{Bryant-1987,Hitchin:2000jd}. This is a nowhere-vanishing $3$-form that
defines a Riemannian metric $g_{\varphi},$ such that for any vectors $u$ and
$v$, the following holds
\begin{equation}
g_{\varphi}\left(  u,v\right)  \mathrm{vol}_{\varphi}=\frac{1}{6}\left(
u\lrcorner\varphi\right)  \wedge\left(  v\lrcorner\varphi\right)
\wedge\varphi. \label{gphi}%
\end{equation}
At any point, the stabilizer of $g_{\varphi}$ (along with orientation) is
$SO\left(  7\right)  $, whereas the stabilizer of $\varphi$ is $G_{2}\subset
SO\left(  7\right)  $. This shows that at a point, positive $3$-forms forms
that correspond to the same metric, i.e., are \emph{isometric}, are
parametrized by $SO\left(  7\right)  /G_{2}\cong\mathbb{RP}^{7}\cong
S^{7}/\mathbb{Z}_{2}$. Therefore, on a Riemannian manifold, metric-compatible
$G_{2}$-structures are parametrized by sections of an $\mathbb{RP}^{7}%
$-bundle, or alternatively, by sections of an $S^{7}$-bundle, with antipodal
points identified. This is precisely the parametrization given by
(\ref{phisameg1}).

Alternatively, a $G_{2}$-structure in a fixed metric class can be interpreted
as a reduction of the principal $SO\left(  7\right)  $-bundle $P$ of
orthonormal frames to a principal $G_{2}$-subbundle, and hence each such
reduction corresponds to a section $\sigma$ of an $SO\left(  7\right)  /G_{2}%
$-bundle $N$ and equivalently, an $SO\left(  7\right)  $-equivariant map
$s:P\longrightarrow SO\left(  7\right)  /G_{2}\cong S^{7}/\mathbb{Z}_{2}$. This
is the picture used in \cite{SaEarpLoubeau}.

We may also use the $G_{2}$-structure $\varphi$ and the metric to define the
\emph{octonion bundle }$\mathbb{O}M\cong\Lambda^{0}\oplus TM$ on $M$ as a rank
$8$ real vector bundle equipped with an octonion product of sections given by
\begin{equation}
A\circ_{\varphi}B=\left(  ab-g\left(  \alpha,\beta\right)  ,a\beta
+b\alpha+\alpha\times_{\varphi}\beta\right)
\end{equation}
for any sections $A=\left(  a,\alpha\right)  $ and $B=\left(  b,\beta\right)
$. We set the metric $g=g_{\varphi}$, since we are fixing the metric, even though the $G_{2}$-structure may change. Here we define $\times_{\varphi}$ by $g\left(  \alpha\times_{\varphi}%
\beta,\gamma\right)  =\varphi\left(  \alpha,\beta,\gamma\right)  $ and given
$A\in\Gamma\left(  \mathbb{O}M\right)  $, we write $A=\left(
\mathop{\rm Re}\nolimits A,\mathop{\rm Im}\nolimits A\right)  .$ The metric on
$TM$ is extended to $\mathbb{O}M$ to give the octonion inner product
$\left\langle A,B\right\rangle =ab+g\left(  \alpha,\beta\right)  $, which is
Hermitian with respect to the octonion product. In the formula
(\ref{phisameg1}), the pair $\left(  a,\alpha\right)  $ can now be interpreted
as a unit octonion section.

The \emph{intrinsic torsion }of a $G_{2}$-structure is defined by
$\nabla\varphi$, where $\nabla$ is the Levi-Civita connection for the metric
$g$ that is defined by $\varphi$. Following \cite{karigiannis-2007}, we have
\begin{equation}
\nabla_{a}\varphi_{bcd}=2T_{a}^{\ e}\psi_{ebcd}^{{}}\ \text{and }\nabla
_{a}\psi_{bcde}=-8T_{a[b}\varphi_{cde]} \label{codiffphi}%
\end{equation}
where $T_{ab}$ is the \emph{full torsion tensor}, note that an additional
factor of $2$ is for convenience, and $\psi=\ast\varphi$ is the $4$-form that
is the Hodge dual of $\varphi$ with respect to the metric $g$. The $G_{2}%
$-structure is known as \emph{torsion-free} if $T=0$, and in that case
$\nabla$ has holonomy contained in $G_{2}$. Conversely, if $\nabla$ has
holonomy contained in $G_{2}$, then there exists a torsion-free $G_{2}%
$-structure within the metric class. Let $V=(a,\alpha)$ be a unit octonion section, then define $\sigma_{V}\left(  \varphi\right)=\sigma_{(a,\alpha)}(\varphi)$, as in (\ref{phisameg1}).  It has been shown in
\cite{GrigorianOctobundle} that the torsion of the $G_{2}$-structure
$\varphi_{V}=\sigma_{V}\left(  \varphi\right)  $ is given by
\begin{equation}
T^{\left(  V\right)  }=VTV^{-1}-\left(  \nabla V\right)  V^{-1}
\label{TV}%
\end{equation}
where $T$ is the torsion of $\varphi,$ interpreted as a $1$-form with values
in the bundle of imaginary octonions $\operatorname{Im}\mathbb{O}M$. If we now
define an octonion covariant derivative $D$ on sections of $\mathbb{O}M$ via
\begin{equation}
DV=\nabla V-VT, \label{DV}%
\end{equation}
the expression (\ref{TV}) simply becomes
\begin{equation}
T^{\left(  V\right)  }=-\left(  DV\right)  V^{-1}. \label{TV2}%
\end{equation}
As shown in \cite{GrigorianOctobundle}$,$ the derivative $D$ has other nice
properties - it is metric-compatible, and satisfies a partial product rule
with respect to octonion product on $\mathbb{O}M$, that is, $D\left(
UV\right)  =\left(  \nabla U\right)  V+U\left(  DV\right)$. Now given
(\ref{TV2}), the divergence of $T^{\left(  V\right)  }$ can be expressed as
\begin{equation}
\operatorname{div}T^{\left(  V\right)  }=-\left(  \Delta_{D}V\right)
V^{-1}-\left\vert DV\right\vert ^{2}. \label{divTV}%
\end{equation}

\section{Energy functional}

\label{secEnergy}Given that the torsion varies across
$G_{2}$-structures within the same metric class, an obvious question is how to
pick a representative of the class with the \textquotedblleft
best\textquotedblright\ torsion. A reasonable way to try and characterize the
best torsion is to look for critical points of a functional. Therefore, given
the set $\mathcal{F}_{g}$ of all $G_{2}$-structures that are compatible with a
given metric $g$, and assuming $M$ is compact, define the functional
$\mathcal{E}:\mathcal{F}_{g}\longrightarrow\mathbb{R}$ by
\begin{equation}
\mathcal{E}\left(  \varphi\right)  =\int_{M}\left\vert T^{\left(
\varphi\right)  }\right\vert ^{2}\operatorname{vol}, \label{Func1}%
\end{equation}
where $T^{\left(  \varphi\right)  }$ is the torsion of a $G_{2}$-structure
$\varphi$. This is the functional used by Dwivedi, Gianniotis, and Karigiannis
in \cite{DGKisoflow}.

As we have seen in the previous section, given a $G_{2}$-structure $\varphi$,
any other $G_{2}$-structure within the same metric class is given by
$\sigma_{V}\left(  \varphi\right)  $ for a unit octonion section $V$.
Therefore, the functional (\ref{Func1}) is equivalent to the functional
$\mathcal{E}_{\mathbb{O}}:\Gamma\left(  S\mathbb{O}M\right)  \longrightarrow
\mathbb{R}$ given by
\begin{equation}
\mathcal{E}_{\mathbb{O}}\left(  V\right)  =\int_{M}\left\vert T^{\left(
V\right)  }\right\vert ^{2}\operatorname{vol}=\int_{M}\left\vert DV\right\vert
^{2}\operatorname{vol} \label{Func2}%
\end{equation}
where we have also applied (\ref{TV2}). Hence, in fact, the functional
$\mathcal{E}_{\varphi}$ is equivalent to an energy functional with respect to
the derivative $D$. This is the functional used in
\cite{GrigorianOctobundle,GrigorianIsoflow}.

On the other hand, following the approach in \cite{SaEarpLoubeau}, recall that
a principal $H$-subbundle of a principal $G$-bundle $P$ may be characterized
by an equivariant map $s:P\longrightarrow$ $G/H$, or equivalently, as a
section $\sigma$ of the associated bundle $N=P\times_{G}\left(  G/H\right)
\cong P/H.$ Assuming that $G$ is semi-simple, so that it admits a bi-invariant
metric, we may define a metric $\eta$ on $N,$ together with the corresponding
Levi-Civita connection $\nabla^{\eta}.$ Moreover, given a metric $g$ on the
base manifold, we may induce a metric on $T^{\ast}M\otimes\sigma^{\ast}TN,$
which is compatible with the splitting $TN=\mathcal{V}N\oplus\mathcal{H}N$
induced by $\nabla^{\eta}.$ Using this metric, we may then define an energy
functional $\mathcal{E}_{\Gamma}:\Gamma\left(  N\right)  \longrightarrow
\mathbb{R}$ on sections of $N$:%

\begin{equation}
\mathcal{E}_{\Gamma}\left(  \sigma\right)  =\int_{M}\left\vert d\sigma
\right\vert ^{2}\operatorname{vol}. \label{Func3}%
\end{equation}
Alternatively, suppose that moreover $G$ is compact, so that $P$ is compact.
Then, let us define an energy functional on $G$-equivariant maps
$s:P\longrightarrow G/H$:%
\begin{equation}
\mathcal{E}_{G}\left(  s\right)  =\int_{P}\left\vert ds\right\vert
^{2}\operatorname{vol}_{P} \label{Func3P}%
\end{equation}
where an induced metric on $T^{\ast}P\otimes s^{\ast}T\left(  G/H\right)  $ is
used. It is then shown in \cite{SaEarpLoubeau}$,$ that for any section
$\sigma\in\Gamma\left(  N\right)  $ and its corresponding $G$-equivariant map
$s\in C_{G}^{\infty}\left(  P,G/H\right)  ,$ $\mathcal{E}_{G}\left(  s\right)
=c_{1}\mathcal{E}_{\Gamma}\left(  \sigma\right)  +c_{2}$ where $c_{1}$ and
$c_{2}$ are uniform constants.

Consider the orthogonal splitting $d\sigma=d^{\mathcal{V}}\sigma
+d^{\mathcal{H}}\sigma$ into horizontal and vertical parts. Since the
horizontal component of the metric is given by $\pi^{\ast}g$, where
$\pi:N\longrightarrow M$ is the bundle projection map, we find that for any
$X\in TM$,
\[
 \left\vert d^{\mathcal{H}}\sigma\left(  X\right)  \right\vert
^{2}=\left(  \pi^{\ast}g\right)  \left(  d\sigma\left(  X\right)
,d\sigma\left(  X\right)  \right)  =g\left(  \left(  \pi\circ\sigma\right)
_{\ast}X,\left(  \pi\circ\sigma\right)  _{\ast}X\right)  =g\left(  X,X\right)
.
\]

Thus, the horizontal part of $d\sigma$ contributes only a constant term to
(\ref{Func3}), and it is thus sufficient to consider just the vertical
component
\begin{equation}
\mathcal{E}_{\Gamma}^{\mathcal{V}}\left(  \sigma\right)  =\int_{M}\left\vert
d^{\mathcal{V}}\sigma\right\vert ^{2}\operatorname{vol}. \label{Func3a}%
\end{equation}
In the $G_{2}$ case, Loubeau and S\'{a} Earp show in \cite{SaEarpLoubeau} that
this functional is equivalent to (\ref{Func1}).

\begin{theorem}
[\cite{SaEarpLoubeau}]If $M$ is $7$-dimensional, $P$ is the $SO\left(  7\right)  $-principal bundle of
oriented orthonormal frames, and $N$ is an associated $SO\left(
7\right)/G_{2}  $-bundle over $M$, then $\left\vert d^{\mathcal{V}}\sigma\right\vert
^{2}=\frac{8}{3}\left\vert T^{\left(  \sigma\right)  }\right\vert ^{2}$ where
$T^{\left(  \sigma\right)  }$ is the torsion tensor of the $G_{2}$-structure
defined by the section $\sigma$.
\end{theorem}

\section{Gradient flow}

\label{secGradient}Given the functionals defined in
the previous section, we may consider critical points and negative gradient
flows of the functionals. This is summarized below.
\[%
\begin{tabular}
[c]{llll}%
\textbf{Space} & \textbf{Functional} & \textbf{Critical points} &
\textbf{Negative gradient flow}\\
$\mathcal{F}_{g}$ & $\mathcal{E}\left(  \varphi\right)  $ &
$\operatorname{div}T^{\left(  \varphi\right)  }=0$ & $\frac{\partial
\varphi_{t}}{\partial t}=2\operatorname{div}T^{\left(  \varphi_{t}\right)
}\lrcorner\psi_{t}$\\
$\Gamma\left(  S\mathbb{O}M\right)  $ & $\mathcal{E}_{\mathbb{O}}\left(
V\right)  $ & $\Delta_{D}V+\left\vert DV\right\vert ^{2}V=0$ & $\frac{\partial
V_{t}}{\partial t}=\Delta_{D}V_{t}+\left\vert DV_{t}\right\vert ^{2}V_{t}$\\
$\Gamma\left(  N\right)  $ & $\mathcal{E}_{\Gamma}\left(  \sigma\right)  $ &
$\tau^{\mathcal{V}}\left(  \sigma\right)  =0$ & $\frac{\partial\sigma_{t}%
}{\partial t}=\tau^{\mathcal{V}}\left(  \sigma_{t}\right)  $\\
$C_{G}^{\infty}\left(  P,G/H\right)  $ & $\mathcal{E}_{G}\left(  s\right)  $ &
$\tau^{\mathcal{H}}\left(  s\right)  =0$ & $\frac{\partial s_{t}}{\partial
t}=\tau^{\mathcal{H}}\left(  s_{t}\right)  $%
\end{tabular}
\ \ \
\]
where $\tau^{\mathcal{V}}\left(  \sigma\right)  :=\operatorname{Tr}_{g}\left(
\nabla^{\eta}d^{\mathcal{V}}\sigma\right)  $ is the \emph{vertical tension}
\emph{field} of the functional $\mathcal{E}_{\Gamma}\left(  \sigma\right)  $
and $\tau^{\mathcal{H}}\left(  s\right)  :=\operatorname{Tr}_{g}^{\mathcal{H}%
}\left(  \nabla^{\eta}ds\right)  $ is the \emph{horizontal} \emph{tension
field }of the functional $\mathcal{E}_{G}\left(  s\right)  .$ It is proved in
\cite[Theorem 1]{ChrisWood1} that $\sigma\in\Gamma\left(  N\right)  $ is a
harmonic \emph{section}, i.e. a critical point of the functional
(\ref{Func3}), if and only if the corresponding $G$-equivariant map $s\in
C_{G}^{\infty}\left(  P,G/H\right)  $ is a \emph{horizontally} \emph{harmonic
map}$,$ that is $\tau^{\mathcal{H}}\left(  s\right)  =0$. In the expression
for $\tau^{\mathcal{H}}\left(  s\right)  $, the trace is just over the
horizontal distribution in $TP.$ It should be emphasized that the reason that
the critical points of $\mathcal{E}_{G}$ are not exactly harmonic maps is that
we are varying over only the equivariant maps, rather than arbitrary maps. On
the other hand, Wood does prove in \cite[Theorem 3]{ChrisWood1}, that if $G/H$
is a normal $G$-homogeneous manifold and the metric on $P\ $is constructed
from any compatible metric on $G$, then $\sigma$ is a harmonic section if and
only $s$ is a harmonic map, that is, $\tau\left(  s\right)
:=\operatorname{Tr}_{g}\left(  \nabla^{\eta}ds\right)  =0.$ Crucially, these
conditions are satisfied for $G=SO\left(  7\right)  $, $H=G_{2}$, and $P$ the
orthonormal frame bundle on $M$. Moreover, as shown in \cite{SaEarpLoubeau},
given these conditions, a family $\sigma_{t}\in\Gamma\left(  N\right)  $
satisfies the harmonic section flow $\frac{\partial\sigma_{t}}{\partial
t}=\tau^{\mathcal{V}}\left(  \sigma_{t}\right)  $ if and only if there is a
corresponding family $s_{t}\in C_{G}^{\infty}\left(  P,G/H\right)  $ that
satisfies the harmonic map flow $\frac{\partial s_{t}}{\partial t}=\tau\left(
s_{t}\right)  .$ Also, Wood has shown in \cite{Wood90} that equivariance
is preserved along the harmonic map flow, so that if the initial condition is
equivariant, then the flow will continue to be equivariant. This shows a close
relationship between harmonic map theory and the theory of harmonic sections,
and hence the flow (\ref{octoflow2}) of isometric $G_{2}$-structures.

On the other hand, one must be careful when applying harmonic map results. In
particular, the energy $\mathcal{E}_{G}\left(  s\right)  $ contains a
topological term that can never be arbitrarily small, and thus standard small
initial energy long-time existence results \cite{ChenDingHM} for harmonic maps
cannot be applied. Similarly, while a constant map is always harmonic, an
equivariant map $s:P\longrightarrow G/H$ can never be constant (if $H \neq G$). Thus existence
of non-trivial harmonic equivariant maps and hence harmonic sections is not
guaranteed, as expected.

Some results from the theory of harmonic maps do carry over, at least in the
$G_{2}$-case. It was shown in \cite{DGKisoflow,GrigorianIsoflow} that almost
monotonicity and $\varepsilon$-regularity results similar to the harmonic map
heat flow \cite{ChenDingHM,StruweChen,StruweHM1} hold for the flow
(\ref{octoflow2}).

Let $p_{x_{0},t_{0}}\left(  x,t\right)  $ be the backward heat kernel on $M$,
that is, the solution of the backward heat equation for $0\leq t\leq t_{0}$
that converges to a delta function at $\left(  x,t\right)  =\left(
x_{0},t_{0}\right)  $. Then, given a time-dependent octonion section $V_{t}$
or equivalently, a $3$-form $\varphi_{t}=\sigma_{V\left(  t\right)  }\left(
\varphi\right)  \ $for some fixed $G_{2}$-structure $\varphi$, define the
$\mathcal{F}$-functional \cite{GrigorianIsoflow}%
\begin{equation}
\mathcal{F}\left(  x_{0},t_{0},t\right)  =\left(  t_{0}-t\right)  \int
_{M}\left\vert T^{\left(  V_{t}\right)  }\left(  x\right)  \right\vert
^{2}p_{x_{0},t_{0}}\left(  x,t\right)  \operatorname{vol}\left(  x\right)  ,
\label{Ffunctional}%
\end{equation}
where $T^{\left(  V_{t}\right)  }=-\left(  DV_{t}\right)  V_{t}^{-1}$ is the
torsion of the $G_{2}$-structure $\varphi_{t}.$ In \cite{DGKisoflow}, the
analogous quantity is denoted by $\Theta_{\left(  x_{0},t_{0}\right)  }\left(
\varphi\left(  t\right)  \right)  $. It is then shown in both \cite[Theorem
5.3]{DGKisoflow} and \cite[Proof of Corollary 7.2]{GrigorianIsoflow} that
$\mathcal{F}$ satisfies an \emph{almost} monotonicity formula along the flow
(\ref{octoflow}). Suppose $V_{t}$ is a solution of the flow (\ref{octoflow})
for $0\leq t<t_{0}$ with initial energy $\mathcal{E}\left(  0\right)
=\mathcal{E}_{0}$. Then, there exists a constant $C>0$, that only depends on
the background geometry, such that for any $t\ $and $\tau$ satisfying
$t_{0}-1\leq\tau\leq t<t_{0}$, $\mathcal{F} $ satisfies the following
relation%
\begin{equation}
\mathcal{F}\left(  x_{0},t_{0},t\right)  \leq C\mathcal{F}\left(  x_{0}%
,t_{0},\tau\right)  +C\left(  t-\tau\right)  \left(  \mathcal{E}%
_{0}+\mathcal{E}_{0}^{\frac{1}{2}}\right)  . \label{Ztmonotonicity}%
\end{equation}
In \cite{DGKisoflow}$,$ the last term in (\ref{Ztmonotonicity}) was $C\left(
t-\tau\right)  \left(  \mathcal{E}_{0}+1\right)  $, which of course follows
from (\ref{Ztmonotonicity}) for a different constant $C$. In both
\cite{DGKisoflow} and \cite{GrigorianIsoflow} similar
versions of an $\varepsilon$-regularity result is proven for $\mathcal{F}$.
We'll state it as in \cite{GrigorianIsoflow}.

\begin{theorem}
[{\cite[Theorem 5.7]{DGKisoflow} and \cite[Theorem 7.1]{GrigorianIsoflow}}%
]\label{thmEpsRegular}Given $\mathcal{E}_{0}$, there exist $\varepsilon
>0\ $and $\beta>0,\ $both depending on $M$ and $\beta$ also depending on
$\mathcal{E}_{0}$, such that if $V$ is a solution of the flow (\ref{octoflow})
on $M\times\lbrack0,t_{0})\ $with energy bounded by $\mathcal{E}_{0}$, and if
\begin{equation}
\mathcal{F}\left(  x_{0},t_{0},t\right)  \leq\varepsilon\label{epsreghypo}%
\end{equation}
for $t\in\lbrack t_{0}-\beta,t_{0}),$ then $V$ extends smoothly to $U_{x_{0}%
}\times\lbrack0,t_{0}]$ for some neighborhood $U_{x_{0}}$ of $x_{0}$ with
$\left\vert DV\right\vert =\left\vert T^{\left(  V\right)  }\right\vert $
bounded uniformly.
\end{theorem}

Then, Theorem \ref{thmEpsRegular} was used in
\cite{DGKisoflow,GrigorianIsoflow} to show long-time existence of the
isometric heat flow and convergence to a $G_{2}$-structure with
$\operatorname{div}T=0$ given sufficiently small initial pointwise torsion.

Given a $G_{2}$-structure $3$-form $\varphi,$ in \cite{DGKisoflow} a concept
of \emph{entropy} was defined:%
\begin{equation}
\lambda\left(  \varphi,\sigma\right)  =\max_{\left(  x,t\right)  \in
M\times\left(  0,\sigma\right)  }\left\{  t\int_{M}\left\vert T^{\left(
\varphi\right)  }\left(  y\right)  \right\vert ^{2}p_{\left(  x,t\right)
}\left(  y,0\right)  \operatorname{vol}\left(  y\right)  \right\}  .
\label{entropy}%
\end{equation}
This mirrors similar entropy concepts defined for the mean curvature flow,
Yang-Mills flow, and the harmonic map heat flow, in
\cite{ColdingMinicozzi2012}, \cite{KelleherStreetsYM1}, and
\cite{BolingKelleherStreetsHM1}, respectively. The quantity $\lambda\left(
\varphi,\sigma\right)  $ is shown in \cite{DGKisoflow} to be invariant under
the scaling $\left(  \varphi,\sigma\right)  \mapsto\left(  c^{3}\varphi
,c^{2}\sigma\right)  .$ While the same quantity could be defined for an
octonion section $V,$ if considered as a function of $V,$ $\lambda$ would lose
the scaling property for $V.$ So in this case, using the $3$-form has an
advantage. Overall, one of the key results in \cite{DGKisoflow} is long term
existence and convergence of the flow (\ref{octoflow2}) given sufficiently
small entropy.

\begin{theorem}
[{\cite[Theorem 5.15]{DGKisoflow}}]Let $\varphi_{0}$ be a $G_{2}$-structure on
a compact $7$-manifold $M$. For any $\delta,\sigma>0$, there exists
$\varepsilon>0$, such that if $\lambda\left(  \varphi_{0},\sigma\right)
<\varepsilon$, then the flow (\ref{octoflow2}) with initial condition
$\varphi\left(  0\right)  =\varphi_{0}$ exists for all time and converges
smoothly to a $G_{2}$-structure $\varphi_{\infty}$ that satisfies
$\operatorname{div}T^{\left(  \varphi_{\infty}\right)  }=0$ and $\left\vert
T^{\left(  \varphi_{\infty}\right)  }\right\vert <\delta.$
\end{theorem}

Although good progress has been made on properties of the flows
(\ref{octoflow}) and (\ref{octoflow2}), many questions still remain. For
example, is it possible to prove long-time existence given small initial
energy, rather than entropy or pointwise torsion? If we combine the
equivariant harmonic map approach with the octonion approach, then everything
could be reformulated in terms of equivariant maps from the orthonormal frame
bundle $P$ to $S^{7}$ equipped with the octonion product. It is likely that
the additional algebraic structure could help achieve stronger results.


{\footnotesize
\bibliographystyle{habbrv}

}
\end{document}